# Dynamic Pricing System for Physical Internet Enabled Hyperconnected Less-than-Truckload Freight Logistics Networks

Tiankuo Zhang, Paria Nourmohammadi, Sixtine Guerin, Benoit Montreuil, Alan Erera
Physical Internet Center, Supply Chain & Logistics Institute,
School of Industrial & Systems Engineering, Georgia Institute of Technology, Atlanta, U.S.A.
Corresponding author: tzhang417@gatech.edu

**Abstract:** *Less-than-truckload (LTL) shipping plays a critical role in modern supply chains by consolidating freight from multiple shippers into shared vehicles. Despite its operational flexibility and potential sustainability benefits, the LTL sector faces persistent challenges, including high per-unit costs and financial instability, as evidenced by recent industry bankruptcies. This paper investigates two structural issues limiting LTL performance: the constrained consolidation potential imposed by proprietary logistics networks, and the inefficiency of fixed pricing models that fail to reflect real-time network conditions. To address these, we explore a Physical Internet (PI)-enabled, hyperconnected LTL logistics system based on open asset sharing and dynamic flow consolidation. We then propose a dynamic pricing framework tailored for this network. Through a simulation-based study grounded in Freight Analysis Framework data and cost estimates from industry sources, we evaluate system performance across three demand and cost uncertainty scenarios in the Southeastern U.S. The results validate our system's effectiveness and suggest a promising path forward for building more efficient LTL logistics operations.*

**Physical Internet (PI) Roadmap Fitness**: ☒ PI Networks, ☒ System of Logistics Networks

**Targeted Delivery Mode-s**: ☒ Paper, ☒ In-Person presentation

## 1 Introduction

Less-than-truckload (LTL) shipping strategies involve consolidating multiple shippers' goods into shared trucks or trailers en route to common intermediate or final destinations. Compared to full truckload (TL) shipments, where a single shipper's freight fills an entire vehicle, LTL provides a flexible alternative that improves truck utilization, reduces empty miles, and supports more sustainable and efficient supply chain operations. As e-commerce and decentralized supply chains continue to grow, the role of LTL in enabling cost-effective, frequent, and distributed shipping has become increasingly central.

However, despite its strategic importance, the LTL industry has faced mounting operational and financial pressures. Recent high-profile bankruptcies, including Yellow Corp's collapse in 2023 and following similar cases in 2019, have underscored the instability within the LTL market. These incidents have cast doubt on LTL carriers' ability to remain profitable in an increasingly competitive and volatile logistics environment. A critical pain point is the cost structure: while LTL offers improved capacity utilization, the cost per pound is significantly higher than that of TL shipments. This disparity is particularly stark when shipments are compared based on similar promised delivery windows, suggesting that LTL customers often pay a premium for similar services.

In this paper, we focus on two key reasons for this phenomenon. The first is that consolidation potential of LTL service providers is bounded by their dedicated network of logistic hubs and their dedicated fleet of trucks and trailers, leading them to constantly compromise between delivery velocity and vehicle fill rate. The second reason is that LTL companies rely on fixed pricing systems with minimal adjustments concerning current network status and shipment details. Therefore, the fixed shipment rates are usually above the underlying real-time rate to stay profitable.

We address the first reason by focusing exclusively on Physical Internet (PI) enabled hyperconnected freight logistic systems for LTL services, based on open asset sharing and open flow consolidation. We propose a PI-LTL hyperconnected freight logistics network based on open asset sharing and open flow consolidation. The services of logistic hubs are openly accessible at a price depending on required sorting, transloading, or crossdocking services, and they can consolidate packages encapsulated in modular PI containers according to intermediary or final destination as required. Similarly, the vehicles used need not to be dedicated to a specific LTL service provider and they can load, move, and unload modular containers under the responsibility of multiple LSPs concurrently. Ultimately, most LTL freight moves in fully loaded trucks and spend minimal dwell time at hubs.

We address the second price-related reason by introducing a dynamic pricing system specifically conceived for hyperconnected LTL logistic systems, with two intelligent agents within a hyperconnected LTL logistics platform. The first pricing agent estimates the shipment price by identifying historical and real-time requests with similar characteristics and computing a robust price using a weighted confidence interval method. The second pricing agent receives the actual optimized cost after a shipment is completed, calculates the deviation, and updates the real-time data pool, enabling continuous learning and refinement. Together, these agents allow the system to adapt pricing dynamically based on evolving network conditions, improving both accuracy and responsiveness over time. This pricing mechanism supports more efficient decision-making while maintaining financial viability for carriers in a shared, open logistics network.

Based on the developed system, we report a simulation-based experiment and examine the pricing system's effectiveness in a hyperconnected LTL network in the Southeastern U.S. and in three scenarios. In scenario 1, we consider a smaller portion of the overall commodity flow with no cost uncertainties, meaning the real-time network has full clairvoyance on arc cost information that is the same with the historical data. In scenario 2, we consider the same commodity flow portion with cost uncertainties, such that the real-time arc costs are different from the historical data, and the arc cost information within the historical data also varies. In scenario 3, we consider a larger commodity flow portion with cost uncertainties. The experiment was set on a multi-day horizon with hourly time instances to spawn requests. It utilizes Freight Analysis Framework (FAF) data to generate representative synthetic LTL freight demand scenarios and to simulate shipment requests between hub pairs. It uses Uber Freight and American Transportation Research Institute (ATRI)'s reports to estimate cost parameter. In each scenario, we compare accumulated profits of our proposed dynamic pricing system versus the fixed pricing strategy. The comparison results validate the dynamic pricing model's ability to increase carrier profitability and market acquisitions

The full paper is organized as follows. Section 2 presents the related literature. Section 3 demonstrates our proposed PI-LTL hyperconnected freight logistics network. Section 4 proposes the framework of our dynamic pricing system. Section 5 analyzes the simulation-based experiment. Section 6 summarizes the contributions, limitations, and future work directions.

## 2 Related Literature

This literature review is divided into two main sections: the first explores research on Physical Internet (PI) networks, while the second focuses on dynamic pricing mechanisms, particularly in less-than-truckload (LTL) logistics.

### 2.1 PI Networks

The Physical Internet (PI) represents a transformative paradigm shift in global logistics, first conceptualized by Montreuil (2011) as a response to the "global logistics sustainability grand challenge," centering on several architectural innovations: hyperconnected network, standardized modular containers, universal interconnectivity protocols, and open logistics networks (Ballot et al., 2014).

Based on these core concepts, researchers have studied applications and impact of PI networks. Hakimi et al. (2012) applied PI concepts and designed an open logistics web in Fance. Through their research, they reported notable improvements in economic, environmental, and social efficiency and performance by adopting PI concepts in real-world logistics systems. Li et al. (2022) designed a PI operating system based on hyperconnected relay network that provides truck drivers with more return-to-domicile opportunities. Other researchers focused on design frameworks of PI networks, aiming to provide conclusive guidance for designing PI networks fitting into different scenarios and criteria. Shaikh et al. (2021) introduced in-hub protocols for dynamically generating consolidation sets of modular containers and requests for on-demand transportation services, ensuring reliable pickup and delivery within promised time windows. Grover et al. (2023) proposed a framework on PI-hyperconnected network design framework that integrates key concepts such as tiered network topology, hub interconnectivity, consolidation, and containerization.

### 2.2 Dynamic Pricing Models in Physical Internet Framework

As the PI ecosystem evolves, pricing mechanisms have emerged as a critical research area. The unique nature of PI environments, characterized by dynamic, stochastic transport requests in open logistics hubs, necessitates innovative pricing approaches distinct from traditional models. Less-than-truckload (LTL) requests with varying volumes and destinations continually arrive and remain available only briefly, creating a complex decision environment for carriers (Sarraj et al., 2014a; Xu & Huang, 2013; Pan et al., 2014). Qiao et al. (2016) made a notable contribution by developing a dynamic pricing model for LTL carriers operating in PI-hubs. Their approach employs dynamic programming and auction theory to optimize bid prices in real-time, adapting to the stochastic arrival of requests with different characteristics. At the core of their research is a pricing estimation function that highlights in capturing opportunity costs of serving one request not others. Yet, this research is limited to only one origin hub and one destination hub, which lacks adaptabilities in large-scale PI networks.

Building on this foundation, Qiao et al. (2020) then expanded the research scope with a more comprehensive model addressing both dynamic pricing and request selection. Their enhanced approach incorporates forecasting and uncertainty and considers the multi-leg transport problem. However, while valuable insights can be drawn from these parallel domains, they often focus on capacity control and acceptance/rejection decisions rather than flexible pricing mechanisms. Moreover, the current carrier-based pricing models in PI research remain rooted in competitive market paradigms that optimize individual carrier revenue rather than system-wide efficiency—a limitation that fails to fully utilize PI-based open asset sharing or realize the collaborative potential of the PI vision.

## 3   Proposed PI-LTL Hyperconnected Freight Logistics Network

At the hub level, once a truck arrives, containers are first sorted based on their next destination. The system prioritizes which containers remain on the truck (those destined for the truck's next stop), while others are unloaded and, later, reallocated to outbound trucks heading toward their respective next stops. If two shipments are bound for the same next destination but only one spot remains, priority is given based on shipment urgency, time window, or contractual service tier. This decision-making process is handled by the digital planning agent, which evaluates available capacity and urgency constraints in real time. The hub operations follow a structured flow: arrival, unloading, reconfiguration, preparation, loading, and departure; each optimized for minimal handling and high-speed turnover, in line with PI container standards and automation compatibility.

At the network level, the PI-LTL (Physical Internet–Less-than-Truckload) network leverages hyperconnected and relay-based transportation, with the key design principle that transportation arcs between hubs should remain under 5.5 hours. This duration enables short-haul relays that comply with driver hour-of-service regulations, reduce driver fatigue, and maximize vehicle and driver productivity through more predictable shift scheduling. The short-haul structure also allows for increased frequency of dispatches and improves service reliability by reducing exposure to long-haul disruptions. Such a hyperconnected and relay-based network framework is modular and highly scalable—it can be applied to any given set of hubs. This modularity allows the network to dynamically grow or shrink by adding or removing hubs, as long as connectivity and capacity thresholds are respected. Each transportation arc supports bidirectional movement, facilitating flexible load balancing and repositioning of assets. Hubs are designed to function as universal access points rather than assets of certain carriers, which enables cooperations between different carriers and logistics service providers. This open-access approach allows for greater integration of regional, national, and even international carriers, promoting a more inclusive and efficient freight system.

The platform also reinforces open asset sharing through centralized orchestration of both capacity and pricing. Instead of relying on traditional competition, where carriers independently bid for individual shipments and customers choose from a fragmented pool, the system acts as a neutral coordinator. It matches shipment requests to available assets (trucks, drivers, and hub capacity) regardless of ownership. This centralized matching engine incorporates real-time data and optimization algorithms to allocate resources fairly and efficiently. Carriers share both trucks and hub capacities and are compensated transparently based on their participation in fulfilling shipment requests, using a usage-based or performance-based revenue-sharing model. This collaborative model leads to improved resource utilization, higher asset uptime, and a more resilient logistics ecosystem that can absorb demand surges or infrastructure disruptions with greater agility. By decoupling service provision from asset ownership, the platform lowers the entry barrier for small and mid-sized carriers, encourages broader participation, and reduces redundancy in fleet and infrastructure deployment. Over time, this structure fosters a more balanced, efficient, and environmentally sustainable freight network—capable of dynamically adapting to evolving customer needs, regulatory constraints, and macroeconomic conditions.

Figure 1 illustrates our designed PI-LTL hyperconnected freight logistics network implemented in our test bed of 8 hubs across Georgia and Florida. The green arcs represent relay transportation arcs (each under 5.5 hours), showcasing the interconnectedness and time-constrained feasibility of the network. The subgraphs on the right depict hub-level PI container handling processes and container movements, highlighting how reconfiguration and dispatch are orchestrated in real-time under PI principles.

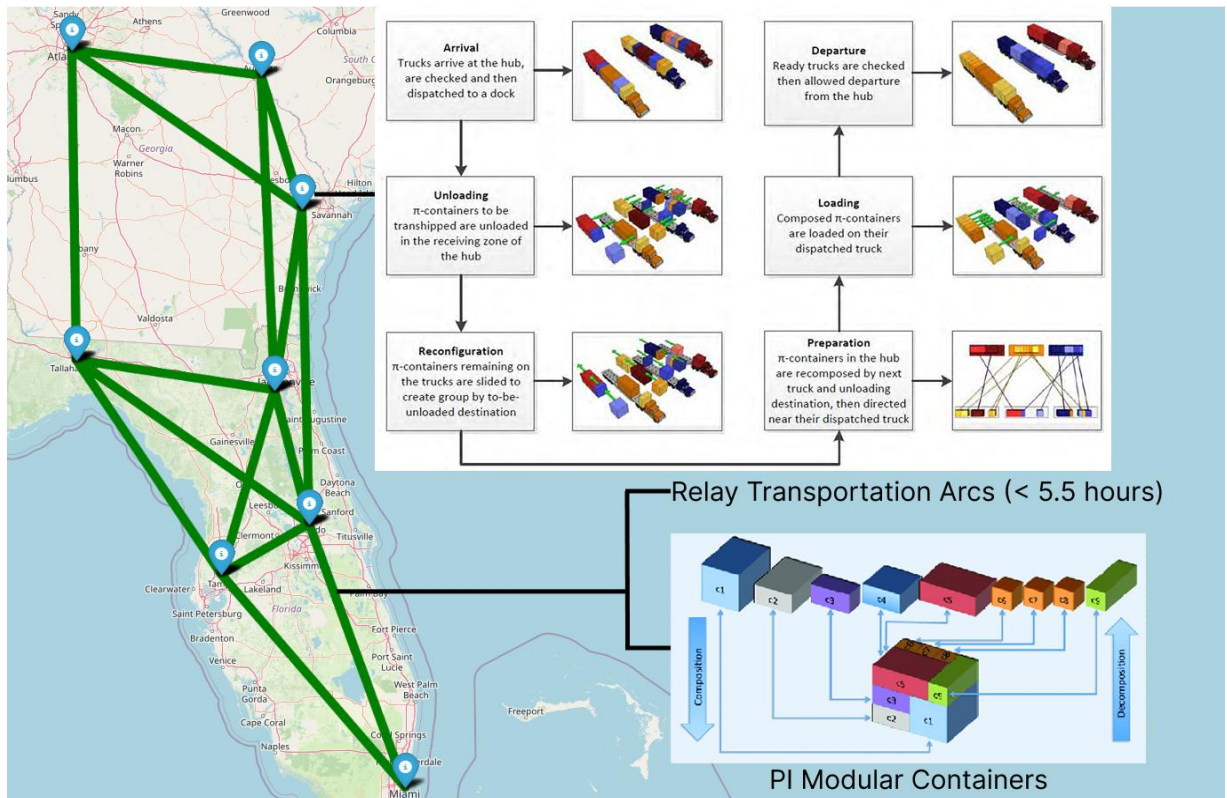

*Figure 1: Illustrating the PI-LTL Hyperconnected Freight Logistics Network*

*(Montreuil, 2011; Montreuil et al., 2013)*

## 4   Proposed Dynamic Pricing System

Our proposed dynamic pricing system is built around two types of intelligent pricing agents. The first pricing agent 1) combines historical data and real-time outcome into shipment request groups based on similar origin, destination, time window, and 2) generates a pricing estimation from the respective request group for the current shipment request based on a given confidence interval for guaranteed robustness. The second type of pricing agent traces the optimized cost of each shipment, calculates the deviation from the absolute difference between the previous estimate and the optimized cost, and logs the outcome to update the real-time outcome pool for the first pricing agent's further use. These two pricing agents are part of a more extensive hyperconnected LTL logistic operating system, as shown in Figure 2, with connections to a consolidation agent and a planning agent that are not the focus of this paper. The above flow chart shows a detailed view following a new shipper request:

1) A shipper asks for a shipment request with specifications to be delivered.
2) The consolidation agent consolidates the request into units of PI modular containers.
3) The first pricing agent finds similar shipments in the system's logs, including historical data and real-time outcomes, and estimates the price of shipping this request.
4) The shipper pays for the estimation, and the planning agent derives the optimal routing and consolidation plan in consideration of open resource sharing and carrier collaboration.
5) After the shipment request is fulfilled, the second pricing agent backtracks the accurate cost of this shipment request and updates logs for further use, and the operating system distribute respective rewards to carriers involved in the shipment requests' fulfillment.

Therefore, as new requests come in and new decisions are made, the pricing system iteratively updates its dynamic pricing estimation matrix with the aims of fast response time, with minimal run-time sacrifice to update results, and increments in accuracy of the pricing estimates.

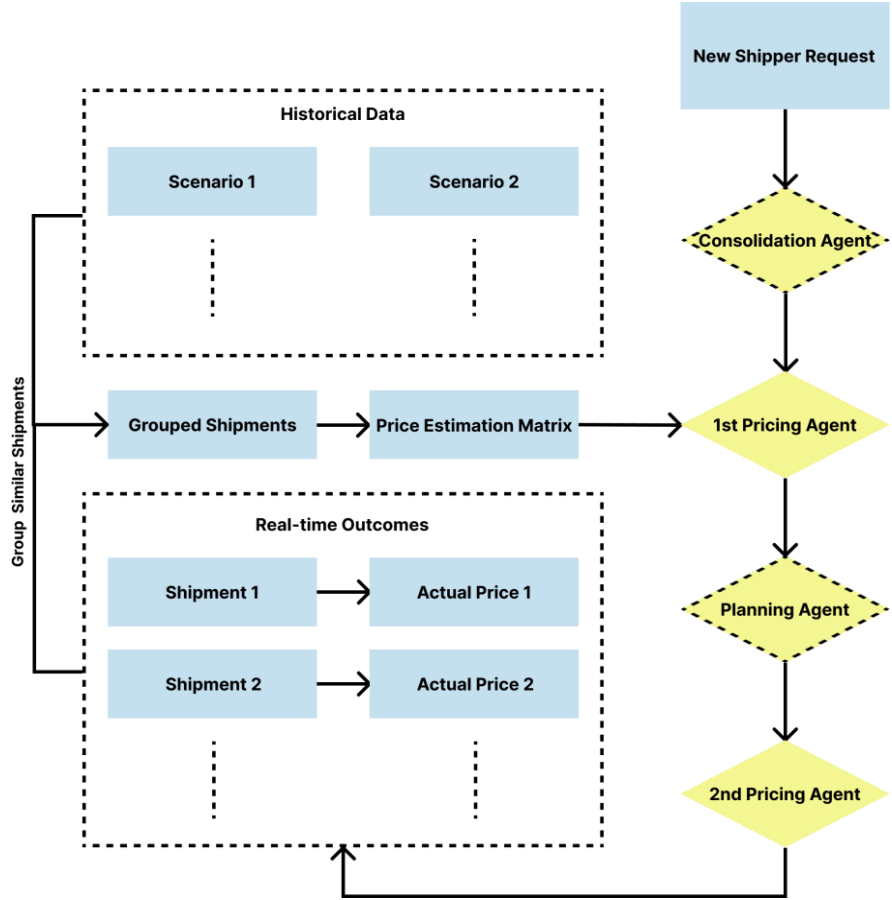

*Figure 2: Demonstrating the Hyperconnected LTL Logistics Operating System*

The first pricing agent estimates price through a weighted confidence interval methodology, in which requests in the historical data and requests in real-time outcomes have different impact to the estimated price. Therefore, in conditions where the current network status is different from the historical data pool, the pricing agents will depend more on information gathered from the current network and can quickly adjust and adapt to the status. For the current shipment request $k$ with origin $o_k$, destination $d_k$, release time $r_k$, deadline $l_k$, and demand $v_k$, the first pricing agent identifies the request's corresponding stored request group $K_s$. For each shipment request $k_s$ in $K_s$, the second pricing agent provides its optimized cost $c_k$ and a confidence weight $w_k$ based on whether $k_s$ belongs to historical data or real-time outcomes. The first pricing agent then calculates the group's weighted mean:

$$\bar{c} = \frac{\sum_{k \in K_s} w_k c_k}{\sum_{k \in K_s} w_k} \tag{1}$$

and the group's weighted variance:

$$s^2 = \frac{\sum_{k \in K_s} w_k (c_k - \bar{c})^2}{\sum_{k \in K_s} w_k} \tag{2}$$

From the weighted mean and weighted variance, the first pricing agent then calculates the group's standard error of the weighted mean as:

$$SE = \sqrt{\frac{s^2}{n_{eff}}} \tag{3}$$

, where $n_{eff} = \frac{\left(\sum_{k \in K_S} w_k\right)^2}{\sum_{k \in K_S} w_k^2}$ is the effective sample size.

With a desired robust level and its corresponding $z$ value, the first pricing agent then gives the robust profitable price $p_k = \bar{c} + z \times SE$ for the current shipment request. Once this request arrives at its destination through the planning agent's optimized route, the second pricing agent will then gather information on its cost $c_k$, which allows us to derive its deviation $\delta_k = |p_k - c_k|$. The second pricing agent then stores this request's origin, destination, time window, demand, and cost information to its real-time outcomes for further use.

We want to point out that the goal of our pricing system is not to have a large revenue $p_k - c_k$ but to guarantee a certain robustness level of remaining profitable while having a low deviation for each request. Therefore, in assessing performances of the pricing system, the desirable outcome is to have low total deviations and low per-pound prices.

## 5 Results and Discussion

We present a simulation-based experiment using a PI-LTL hyperconnected freight logistics network based on freight zones from the FAF dataset, and we focus on 8 zones in Georgia and Florida. Besides 8 hubs derived from the 8 zones, the PI-hyperconnected freight logistics network utilizes relay transportation and includes 30 transportation arcs that are less than 5.5 hours, such that drivers can return to their based hub abiding the driving regulations. We consider a time horizon of 48 hours with hourly planning instances and derive commodity requests with varied assumptions in three scenarios as shown in Table 1. For robustness of our results, we generate 100 cases as historical data and 30 cases as testing cases for each scenario. In this experiment, the target profit robustness is 90% in all cases of each scenario, which can be easily changed to adapt different goals.

*Table 1: Comparison between setups of three scenarios*

| Scenario | Number of Cases | Percentage of LTL Flow | Volume Upper Bound | Per-Truck-Mile Cost | Per-Pound Transload Cost |
|---|---|---|---|---|---|
| 1 | 100/30 | 1% | 1/3 Truckload | $2 | $0.04 |
| 2 | 100/30 | 1% | 1/3 Truckload | $1-$3 | $0.02-$0.06 |
| 3 | 100/30 | 5% | Full Truckload | $1-$3 | $0.02-$0.06 |

In the first scenario, we assume the PI-LTL network handles 1% of all LTL shipment flows (10% of the entire shipment flow), which is approximately 500 to 600 full truckloads in a 48-hour time horizon. From the estimated LTL shipment flows handled by the PI-LTL network, we randomly generate individual requests with less than 1/3 of a full truckload and time windows between the shortest travel time plus 4 hours and one day for creating feasible yet time-sensitive requests. To our best knowledge, this time window interval is tighter than most industrial standards, including large and established logistics providers like UPS and USPS. As we consider no cost fluctuations in the first scenario, in each case, each traveling arc has the same per truck cost, derived from 2023 trucking operational cost analysis of ATRI and the 2024 annual freight trucking rate report of Uber Freight. Besides, each shipment request also has the same handling cost on each traveling arc in each case.

In the second scenario, the PI-LTL network handles the same amount of LTL shipment flows, and requests and their specifications are randomly generated with the same logic. However, we introduce uniform cost fluctuations, so, in each case, the same traveling arc might have different per-truck cost and different handling cost for the same request. This scenario applies our pricing agents in more practical cases where full pricing information is not disclosed, such that, in testing cases, the pricing agent will estimate cost with historical cost information from different network status.

In the third scenario, the PI-LTL network handles 5% of all LTL shipment flows, which is approximately 2500 to 3000 full truckloads in a 48-hour period. The request generated can now take up to full truckload. The cost also fluctuates between each case, following the logic of the second scenario. This scenario resembles cases when the PI-LTL network gains more market share and breaks the traditional concept that LTL shipments only take up to 1/3 of a full truckload, which are resembled by the previous two scenarios. Figure 3 compares the total shipment flows in each scenario, with scenario 1 and 2's shipment flows on the left, and scenario 3's shipment flows on the right.

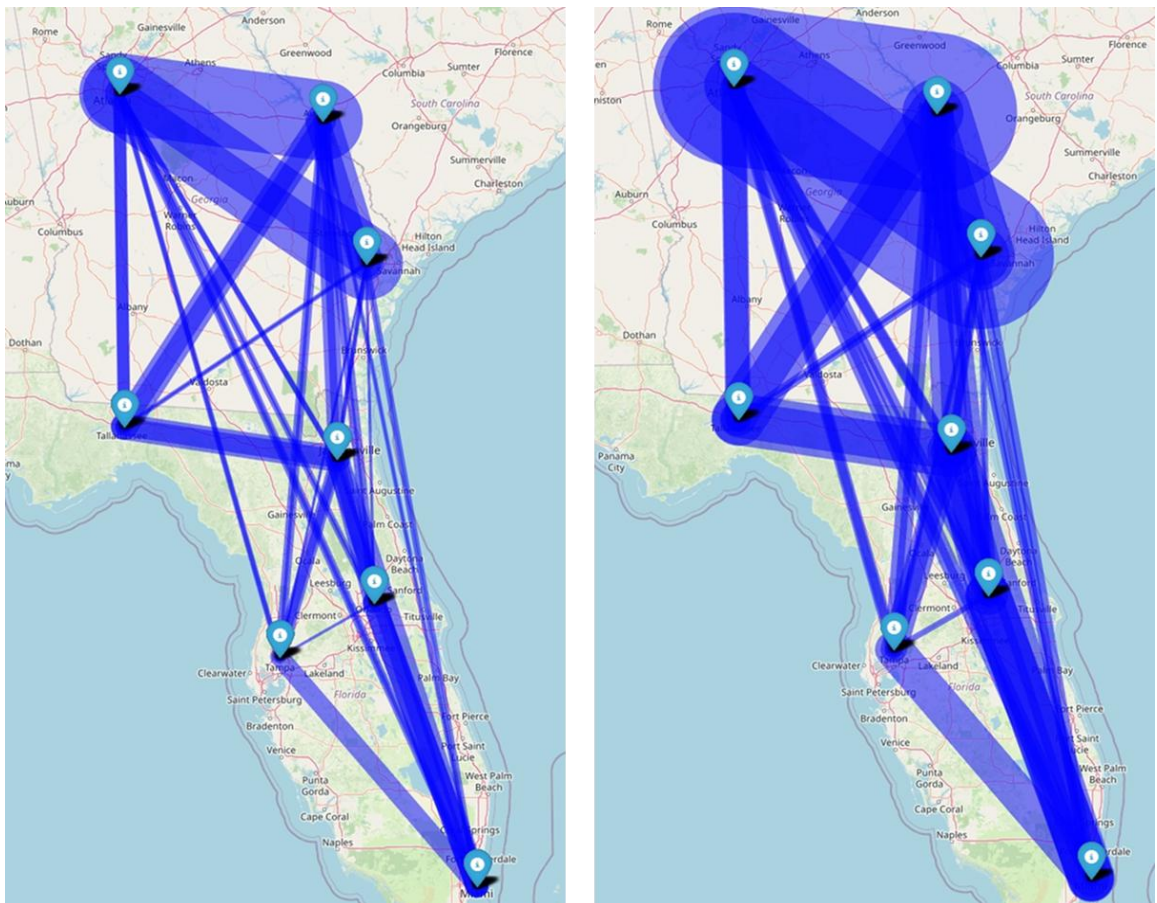

*Figure 3: Comparing Total Shipment Flows Between Scenario 1 (left), 2(left), and 3 (right)*

Table 2 compares the computed results of each scenario with three key performance indicators (KPIs): total deviation, per pound deviation, and average price per pound. Total deviation accounts for each case's sum of deviations of all shipment requests' estimated price from the optimized price; per pound deviation accounts for each case's total deviation divided by the total shipment weight; average price per pound accounts for each case's total price estimation divided by the total weight. As we have multiple cases for each scenario for robustness of our result, we analyze each KPI in lower quartile (25%), median, and upper quartile (75%) at the case level.

*Table 2: Comparison between results of three scenarios*

**Scenario 1**

| KPIs | Lower Quartile | Median | Upper Quartile |
|---|---|---|---|
| Total Deviation | 58320 | 67689 | 74861 |
| Per Pound Deviation | 0.11 | 0.12 | 0.14 |
| Average Price Per Pound | 0.53 | 0.54 | 0.56 |

**Scenario 2**

| KPIs | Lower Quartile | Median | Upper Quartile |
|---|---|---|---|
| Total Deviation | 75601 | 84794 | 109198 |
| Per Pound Deviation | 0.14 | 0.16 | 0.20 |
| Average Price Per Pound | 0.53 | 0.55 | 0.57 |

**Scenario 3**

| KPIs | Lower Quartile | Median | Upper Quartile |
|---|---|---|---|
| Total Deviation | 336353 | 476058 | 593327 |
| Per Pound Deviation | 0.13 | 0.18 | 0.22 |
| Average Price Per Pound | 0.52 | 0.53 | 0.53 |

Comparing scenario 2 to scenario 1, adding uncertainties to the PI-LTL network increases the value of all KPIs, which can be expected as uncertainties will negatively impact our pricing system's prediction. Comparing scenario 3 to the previous two scenarios, due to the increment in request numbers and total volume, the total deviation inevitably increases. However, scenario 3's per-pound average price is below scenario 1 and 2's counterpart. This comparison suggests that our PI-LTL network exhibits an economy of scales. A potential interpretation is that, as our PI-LTL system attracts more shipment requests of larger volumes, the system can have a larger volume of historical knowledge for pricing estimations and explore more consolidation opportunities, the prediction error and prices will be driven down.

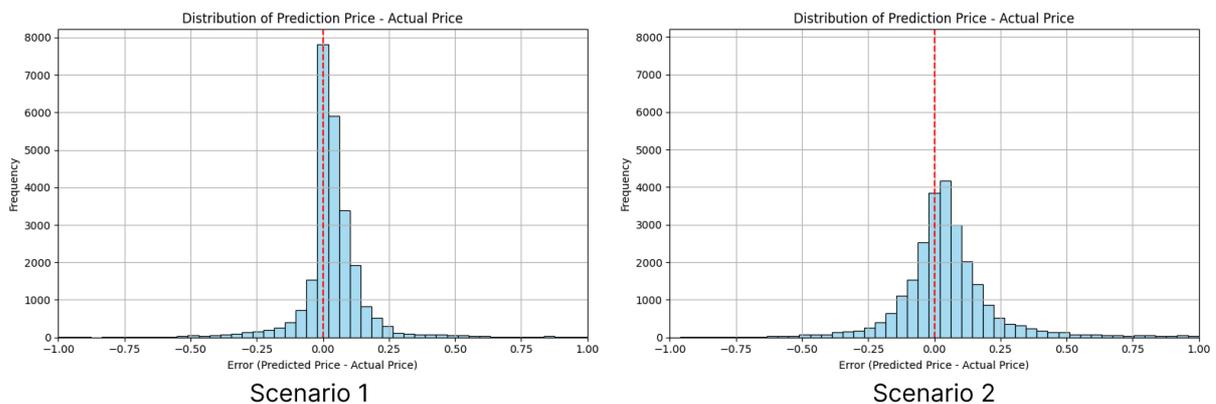

*Figure 4: Histograms of Distributions of Pricing Errors Between Scenario 1 and Scenario 2*

Figure 4 compares the distribution of pricing errors (prediction per-pound price minus the actual per-pound price) between scenario 1 and 2. As scenario 2 introduces uncertainty, less counts of shipment requests' prices are accurately predicted as in scenario 1. However, as shown in scenario 2's histogram, more shipment requests are on the nonnegative half, indicating our system's success in assuring profitability facing uncertainties.

# 6 Conclusion

Our study has three main contributions. First, it proposes a PI-LTL hyperconnected freight logistics system with open asset sharing, smart in-hub consolidation strategies, and hyperconnected relay transportation. Second, it provides a paired dynamic pricing system, as a part of a more advanced LTL operating system, to estimate price based on each shipment request's specifications, historical data, and current real-time outcomes. Third, it conducts a simulation-based experiment in a test bed of Georgia and Florida with real-life industrial data.

This research also opens several avenues for future work. The first direction is to consider a specific profit-splitting strategy to distribute profits fairly among the carriers involved in each shipment request's fulfillment. The second direction is to develop more advanced machine learning models to predict prices with faster response time and higher accuracy. The third direction involves considering a differentiation between premium requests and normal requests with different pricing tiers, such that premium requests are guaranteed to arrive on time with extra cost, and normal requests are less expensive but may have risks of arriving late.

## References

- Montreuil B. (2011): Towards a Physical Internet: Meeting the Global Logistics Sustainability Grand Challenge, Logistics Research, v3, no2-3, 71-87, https://doi.org/10.1007/s12159-011-0045-x
- Montreuil, B., Meller, R. D., & Ballot, E. (2012). Physical Internet foundations. *IFAC Proceedings Volumes*, *45*(6), 26–30. https://doi.org/10.3182/20120523-3-RO-2023.00444
- Hakimi, D., Montreuil, B., Sarraj, R., Ballot, E., & Pan, S. (2012). Simulating a physical internet enabled mobility web: the case of mass distribution in France. In 9th International Conference on Modeling, Optimization & SIMulation-MOSIM'12 (pp. 10-p).
- Li, J., Montreuil, B., & Campos, M. (2022). Trucker-sensitive hyperconnected large-freight transportation: An operating system. In IISE Annual Conference. Proceedings (pp. 1-6). Institute of Industrial and Systems Engineers (IISE).
- Shaikh, S. J., Montreuil, B., Hodjat-Shamami, M., & Gupta, A. (2021). Introducing Services and Protocols for Inter-Hub Transportation in the Physical Internet. arXiv preprint arXiv:2111.07520.
- Grover, N., Shaikh, S. J., Faugère, L., & Montreuil, B. (2023). Surfing the Physical Internet with Hyperconnected Logistics Networks. Athens, Greece. 9th International Physical Internet Conference.
- Qiao, B., Pan, S., & Ballot, E. (2020). *Revenue optimization for less-than-truckload carriers in the Physical Internet: Dynamic pricing and request selection*. Computers & Industrial Engineering, 139, 105563. https://doi.org/10.1016/j.cie.2018.12.010
- Qiao, B., Pan, S., & Ballot, E. (2017). Dynamic pricing model for less-than-truckload carriers in the Physical Internet. *Journal of Intelligent Manufacturing, 30*(6), 2631–2643. https://doi.org/10.1007/s10845-016-1289-8
- Sarraj, R., Ballot, E., Pan, S., & Montreuil, B. (2014). Interconnected logistic networks and protocols: Simulation-based efficiency assessment. *International Journal of Production Research, 52*(11), 3185–3208. https://doi.org/10.1080/00207543.2013.865853
- Montreuil, B., Meller, R. D., Thivierge, C., & Montreuil, Z. (2013). Functional design of physical internet facilities: a unimodal road-based crossdocking hub. Montreal, Canada: CIRRELT.